\documentclass[conference]{IEEEtran}
\IEEEoverridecommandlockouts
\usepackage{cite}
\usepackage[dvipsnames]{xcolor}
\usepackage{algorithmic}
\usepackage{textcomp}
\usepackage{xcolor}
\usepackage{mathpazo}
\usepackage{times}
\usepackage{enumitem}
\usepackage{adjustbox}
\usepackage{amsmath}
\usepackage{amsfonts}
\usepackage{latexsym}
\usepackage{amssymb}
\usepackage{mathabx}
\usepackage{float}
\usepackage{wrapfig}
\usepackage{upref}
\usepackage{theorem}
\usepackage{graphicx}
\usepackage{psfrag}
\usepackage{cite}
\usepackage{ulem}

\usepackage{color}

\usepackage{tikz}
\usetikzlibrary{arrows.meta,bending}

\makeatletter
\newcommand\curvearrowed@[3]
  {%
    \begin{tikzpicture}
      \node[inner sep=.05ex](a){\kern-.25ex#3\kern-.05ex};
      \draw[arrows={-Latex[#1]},line width=#2]
        (a.north west) to[bend left=45] (a.north east);
    \end{tikzpicture}%
  }
\newcommand\curvearrowed[1]
  {%
    \relax
    \ifmmode
      \mathchoice
        {\curvearrowed@{}{.4pt}{$\displaystyle #1$}}
        {\curvearrowed@{}{.4pt}{$\textstyle #1$}}
        {\curvearrowed@{scale=.8}{.325pt}{$\scriptstyle #1$}}
        {\curvearrowed@{scale=.7}{.25pt}{$\scriptscriptstyle #1$}}%
    \else
      \curvearrowed@{}{.4pt}{#1}%
    \fi
  }
\makeatother

\newcommand{\removelatexerror}{\let\@latex@error\@gobble}
\makeatletter
\newcommand{\proofpart}[2]{%
	\par
	\addvspace{\medskipamount}%
	\noindent\emph{Part #1: #2}\par\nobreak
	\addvspace{\smallskipamount}%
	\@afterheading
}
\makeatother

\theoremstyle{plain}
\theorembodyfont{\normalfont\slshape}

\newtheorem{thm}{Theorem$\!$}
\newenvironment{theorem}
{\begin{thm}\hspace*{-1ex}{\bf.}}{\end{thm}}

\newtheorem{clm}[thm]{Claim$\!$}

\newtheorem{lem}[thm]{Lemma$\!$}
\newenvironment{lemma}{\begin{lem}\hspace*{-1ex}{\bf.}}{\end{lem}}

\newtheorem{prop}[thm]{Proposition$\!$}

\newtheorem{cor}[thm]{Corollary$\!$}
\newenvironment{corollary}{\begin{cor}\hspace*{-1ex}{\bf.}}{\end{cor}}

\newtheorem{defn}[thm]{Definition$\!$}

\newtheorem{xmpl}[thm]{Example$\!$}

\newtheorem{cnstr}{Construction$\!$}

\newtheorem{rmk}[thm]{Remark$\!$}

\newcounter{enumrom}
\renewcommand{\theenumrom}{(\roman{enumrom})}

\makeatletter
\renewcommand{\@endtheorem}{\endtrivlist}
\makeatother

\makeatletter
\renewcommand{\thefigure}{{\@arabic\c@figure}}
\renewcommand{\fnum@figure}{{\bf Figure\,\thefigure}}
\makeatother

\newcommand{\cB}{\mathcal{B}}

\newcommand{\cD}{\mathcal{D}}

\newcommand{\cN}{\mathcal{N}}

\newcommand{\bfp}{{\boldsymbol p}}
\newcommand{\bfq}{{\boldsymbol q}}

\newcommand{\bfw}{{\boldsymbol w}}
\newcommand{\bfx}{{\boldsymbol x}}
\newcommand{\bfy}{{\boldsymbol y}}
\newcommand{\bfz}{{\boldsymbol z}}

\renewcommand{\leq}{\leqslant}

\renewcommand{\geq}{\geqslant}

\newcommand{\Cref}[1]{Co\-ro\-lla\-ry\,\ref{#1}}

\usepackage{comment}

\usepackage{amsmath,graphicx}

\outer\def\proclaim #1. #2\par{\medbreak
 \noindent{\bf#1.\enspace}{\sl#2\par}%
 \ifdim\lastskip<\medskipamount \removelastskip\penalty55\medskip\fi}

\begin{document}
\title{\textbf{The Gapped $k$-Deck Problem}
\thanks{The work was funded by NSF grant 2008125, Coded String Reconstruction Problems in Molecular Storage. In the author list, $\ddagger$ denotes equal contribution.}
}
\makeatletter
\newcommand{\linebreakand}{%
  \end{@IEEEauthorhalign}
  \hfill\mbox{}\par
  \mbox{}\hfill\begin{@IEEEauthorhalign}
}
\makeatother
\author{
\IEEEauthorblockN{Rebecca Golm$\ddagger$}
\IEEEauthorblockA{\textit{ECE Department} \\
\textit{University of Illinois}\\
rgolm2@illinois.edu}
\and
\IEEEauthorblockN{Mina Nahvi$\ddagger$}
\IEEEauthorblockA{\textit{Department of Mathematics} \\
\textit{University of Illinois}\\
mnahvi2@illinois.edu} 
\and
\IEEEauthorblockN{Ryan Gabrys}
\IEEEauthorblockA{\textit{Calit2} \\
\textit{University of California, San Diego}\\
ryan.gabrys@gmail.com}
\and
\IEEEauthorblockN{Olgica Milenkovic}
\IEEEauthorblockA{\textit{ECE Department} \\
\textit{University of Illinois}\\
milenkov@illinois.edu}
}
\maketitle
\begin{abstract} The $k$-deck problem is concerned with finding the smallest positive integer $S(k)$ such that there exist at least two strings of length $S(k)$ that share the same $k$-deck, i.e., the multiset of subsequences of length $k$. We introduce the new problem of gapped $k$-deck reconstruction: For a given gap parameter $s$, we seek the smallest positive integer $G_s(k)$ such that there exist at least two distinct strings of length $G_s(k)$ that cannot be distinguished based on a ``gapped'' set of $k$-subsequences. The gap constraint requires the elements in the subsequences to be at least $s$ positions apart within the original string. Our results are as follows. First, we show how to construct sequences sharing the same $2$-gapped $k$-deck using a nontrivial modification of the recursive Morse-Thue string construction procedure. This establishes the first known constructive upper bound on $G_2(k)$. Second, we further improve this bound using the approach by Dudik and Schulman~\cite{Dud}. 
\end{abstract}

\section{Introduction}
\vspace{-0.05in}

The problem of reconstructing strings based on evidence sets of the form of subsequences, substrings or weights of substrings has received significant attention from the theoretical computer science, bioinformatics, and information theory communities alike~\cite{Kal,Batu,Margaritis,Ukkonen,Acharya,Kiah,Chase,Gabrys}. One special instance of this class of problems is the $k$-deck problem~\cite{Kal,scott,Dud,GM,Chrisnata,Chase}, of interest due to its connection to trace reconstruction~\cite{Batu,Cher} and its applications in DNA-based data storage~\cite{Yazdi}.

For a string $\bfx$ of length $n$, the multiset of the $\binom{n}{k}$ subsequences (i.e., ordered collections of not necessarily adjacent entries) of $\bfx$ of length $k$ is called the \textit{$k$-deck} of $\bfx$. We say that $\bfx$ is \textit{$k$-reconstructible} if it is uniquely determined by its $k$-deck, meaning that there exists no other string that has the same $k$-deck as $\bfx$. For example, $(1,0,0,1)$ and $(0,1,1,0)$ have the same $2$-deck, and are hence not $2$-reconstructible. A simple counting argument shows that if two sequences $\bfx$ and $\bfy$ have the same $k$-deck, they also have the same $l$-deck for all $1\leq l\leq k$.

Let $S(k)$ be the smallest positive integer $n$ such that there exist two distinct strings of length $n$ with the same $k$-deck. Kalashnik~\cite{Kal} raised the question of determining $S(k)$. Manvel, Meyerowitz, Schwenk, Smith and Stockmeyer~\cite{Manv} showed that $2k\leq S(k)\leq 2^k$. They proved the upper bound as follows. For two strings $\bfx$ and $\bfy$ of length $n$, let $\bfx\bfy$ be the string obtained by concatenating $\bfx$ and $\bfy$ (note that when concatenating a single bit, say $0$, and a string $\bfx$, we also use the notation $(0,\bfx)$). If $\bfx$ and $\bfy$ have the same $k$-deck, then $\bfx\bfy$ and $\bfy\bfx$ have the same $(k+1)$-deck. The upper bound follows immediately when coupled with the fact that $(0,1)$ and $(1,0)$ have the same $1$-deck. The construction is often referred to as the Morse-Thue construction and the resulting strings are the well-known Morse-Thue strings~\cite{All}. Furthermore, the authors of~\cite{Manv} also showed that in order to prove that every string of length $n$ is $k$-reconstructible, it is enough to prove that every binary string of length $n$ is $k$-reconstructible. Dudik and Schulman~\cite{Dud} improved the above upper bound on $S(k)$ to $\exp{(\frac{3+o(1)}{2\log3}\log^2 k)}$. In the literature, both bounds on the smallest $k$ and $n$ (for a given $n$ and $k$, respectively) for unique and nonunique $k$-deck reconstruction have been reported. 

We define the \textit{gapped $k$-deck} of a binary string $\bfx$ as the multiset of all subsequences of length $\leq k$ that do not include two consecutive entries in $\bfx$. This definition can be extended to larger gaps between entries in $\bfx$: The \textit{$s$-gapped $k$-deck} of a binary string $\bfx$ is the multiset of all subsequences $(x_{i_1},\ldots,x_{i_{\ell}}), 1\leq\ell \leq k,$ such that for all $1\leq j\leq \ell-1$, we have $i_{j+1}\geq i_j+s$. With this definition, the gapped $k$-deck reduces to the $2$-gapped $k$-deck. The problem of interest is to bound $G_s(k)$, the smallest positive integer $n$ for which there exist two binary strings that share the same $s$-gapped $k$-deck. For simplicity, when $s=2$, we write $G(s)$ and refer to the corresponding setting as the \textit{gapped} $k$-deck. Note that unlike the case without gaps, two strings $\bfx$ and $\bfy$ having the same multiset of gapped subsequences of length $k$ does not imply that they also have the same multiset of gapped sequences of length $l$ for some $l<k$. For example, the strings $(0,1,1,1,0)$ and $(1,0,0,0,1)$ have the same multiset of gapped subsequences of length $2$, but they clearly have different multisets of gapped subsequences of length $1$ (which by definition, is the multiset of bits (composition) of the strings). The gapped $k$-deck problem is of interest in molecular storage systems for which readouts are based on nanopore technologies, in which ``gaps'' in readouts arise due to skipping effects~\cite{Yazdi}.

We initiate the study of reconstruction limits of strings given their $s$-gapped $k$-decks and present the first upper bounds on $G_s(k)$ and $G(k)$ in particular. In Section~\ref{sec:preliminaries} we provide necessary preliminaries, while in Section~\ref{sec:mt} we describe a nontrivial extension of a Morse-Thue type construction for $2$-gapped $k$-decks. In Section~\ref{sec:mt}, we state the result for general values of $s$ but omit the proof. Section~\ref{sec:improve} presents an improvement of the upper bound for $G(k)$ from Section~\ref{sec:mt}, based on an adaptation of the method described in~\cite{Dud}.
\vspace{-0.05in} 
\section{Preliminaries} \label{sec:preliminaries}
For a string $\bfx=(x_1, \ldots, x_n) \in \{0,1\}^n$, let 
\begin{align}\label{eq:gkdeck}
\cB^{(k)}(\bfx) = \{{ \left( x_{i_1}, x_{i_2}, \ldots, x_{i_\ell}  \right) : i_j \geq i_{j-1}+2, 0\leq\ell \leq k \}}
\end{align}
denote the multiset of all subsequences of $\bfx$ of length $\leq k$ such that the index of every entry used in a subsequence is nonadjacent in the original string. Also, let
\begin{align*}
\cD^{(k)}(\bfx) = \{{ \left( x_{i_1}, x_{i_2}, \ldots, x_{i_k}  \right) : i_j \geq i_{j-1}+2\}},
\end{align*}
be the \textit{exact} gapped $k$-deck of $\bfx$. Here, we assume that $\cB^{(0)}(\bfx)=\cD^{(0)}(\bfx)=\emptyset$. Clearly, $\cB^{(k)}(\bfx)=\bigcup_{i=0}^k \cD^{(i)}(\bfx).$ As mentioned in the introduction, unlike the classical (ungapped) case, the problem of reconstructing $\bfx$ from $\cB^{(k)}(\bfx)$ differs from that of reconstructing $\bfx$ from $\cD^{(k)}(\bfx)$. Our focus is on finding $G(k)$, the smallest integer $n$ such that there exist two distinct binary strings of length $n$ with the same gapped $i$-deck for all $1\leq i\leq k$. Alternatively, $G(k)$ is the smallest integer $n$ such that there exist two distinct binary strings of length $n$, $\bfx$ and $\bfy$, satisfying $\cB^{(k)}(\bfx)=\cB^{(k)}(\bfy)$.
It is worth pointing out that if $n$ is the smallest integer such that there exist two strings $\bfx$ and $\bfy$ of length $n$ with $\cD^{(k)}(\bfx)=\cD^{(k)}(\bfy)$, then $n=2k-1$. We have $n\geq 2k-1$ because a string of length less than $2k-1$ has no gapped subsequence of length $k$. On the other hand, for any string $\bfz=(z_1\ldots z_k)$ of length $k$, all the strings of length $2k-1$ of the form $(z_1x_1z_2x_2\ldots z_{k-1}x_{k-1}z_k)$ have the same gapped $k$-deck because the only gapped $k$-subsequence of $\bfx$ is $\bfz$. This observation generalizes for $s$-gapped $k$-decks and $n=sk-1$.

In Section~\ref{sec:mt}, we prove that $G(k)\leq 4(2^k-1)-2$. We also provide an upper bound on $G_s(k)$, the smallest integer $n$ such that there exist two distinct strings of length $n$ with the same $s$-gapped $i$-deck for all $1\leq i\leq k$, where $s\geq2$. The bound reads as $G_s(k)\leq (5s-2)2^{k-1}-5s+4$, but the accompanying proof is omitted due to space limitations. The proof of our first bound on $G(k)$ builds upon the next lemma.\vspace{-0.045in}
\begin{lemma}\cite{Manv}\label{manvel}
If $\bfx=(x_1,x_2,\cdots, x_m)$ and $\bfy=(y_1,y_2,\cdots, y_m)$ have the same $k$-deck, then the two concatenation strings $\bfx\bfy=(x_1,\cdots, x_m, y_1,\cdots, y_m)$ and $\bfy\bfx=(y_1,\cdots, y_m,x_1,\cdots, x_m)$ have the same $(k+1)$-deck.
\end{lemma}
\vspace{-0.05in}
\begin{IEEEproof}
The following correspondence proves the claim: Pick any subsequence $\bfz$ of $\bfx\bfy$ of length at most $k+1$. If $\bfz$ is fully contained within the $\bfx$ (or $\bfy$) substring, let $\phi(\bfz)$ be the same subsequence in the $\bfx$ (or $\bfy$) substring of $\bfy\bfx$. Now, assume $\bfz=\bfz_1\bfz_2,$ where $\bfz_1$ is a subsequence of $\bfx$ and $\bfz_2$ is a subsequence of $\bfy$. Note that $\bfz_1$ and $\bfz_2$ have length at most $k$, therefore there exists a subsequence $\bfw_1$ of $\bfy$ that equals $\bfz_1$, due to the fact that $\bfx$ and $\bfy$ have the same $i$-deck for all $1\leq i\leq k$. Similarly, $\bfx$ contains a subsequence $\bfw_2$ that equals $\bfz_2$. Now, let $\phi(\bfz)=\bfw_1\bfw_2$. Therefore, $\bfx\bfy$ and $\bfy\bfx$ have the same $(k+1)$-deck.
\end{IEEEproof}
Using the strings $\bfx=(0,1)$, $\bfy=(1,0)$ and $k=1$ to initialize the recursion, we can see that $(0,1,1,0)$ and $(1,0,0,1)$ have the same $2$-deck. Repeating the process, we find that $(0,1,1,0,1,0,0,1)$ and $(1,0,0,1,0,1,1,0)$ have the same $3$-deck and so on. However, this construction does not work for the gapped case. For example, $(1,0)$ and $(0,1)$ have the same gapped $1$-deck (i.e., composition), but $(1,0,0,1)$ and $(0,1,1,0)$ do not have the same gapped $2$-deck. The reason why the construction fails is that we cannot pick both $x_m$ and $y_1$ (as defined in Lemma~\ref{manvel}) when choosing a gapped subsequence of $\bfx\bfy$. Hence, we need to ``pad'' the boundary between the two concatenated strings in an adequate manner.

\section{The Padded Morse-Thue Sequence Approach} \label{sec:mt}
We prove the existence of two strings $\bfx, \bfy \in \{0,1\}^n$, where $n=4(2^k - 1)-2$, that satisfy $\cB^{(k)}(\bfx) = \cB^{(k)}(\bfy)$, using induction. We start with a few definitions. For a binary string $\bfx=(x_1, \ldots, x_n) \in \{0,1\}^n$ let
\begin{align}
\cB^{(k)}_{L}(\bfx) &:= \cB^{(k)} \left( x_2, x_3, \ldots, x_n \right), \label{eq:defs1}\\
\cB^{(k)}_{R}(\bfx) &:= \cB^{(k)} \left( x_1, x_2, \ldots, x_{n-1} \right),\label{eq:defs2}\\
\cB^{(k)}_{LR}(\bfx) &:= \cB^{(k)} \left( x_2, \ldots, x_{n-1} \right). \label{eq:defs3}
\end{align}
Note that (\ref{eq:defs1}) represents the multiset of all gapped subsequences formed by puncturing $\bfx$ on the left, (\ref{eq:defs2}) represents the multiset of all gapped subsequences formed by puncturing $\bfx$ on the right, while (\ref{eq:defs3}) represents the multiset of all gapped subsequences formed by puncturing $\bfx$ on both ends. We define the sets $\cD^{(k)}_{L}(\bfx), \cD^{(k)}_{R}(\bfx)$, and $\cD^{(k)}_{LR}(\bfx)$ analogously.

We initialize two strings for the ``degenerate" case of $k=1$, corresponding to equal compositions, as follows:
\begin{align}
\bfx^{(1)} &= \left(0,0,1,0 \right ), \hspace{3mm} \bfy^{(1)} =  \left(0,1,0,0 \right )\label{eq:x1}.
\end{align}
Puncturing the first bit from both $\bfx^{(1)}$ and $\bfy^{(1)}$ produces strings that still share the same gapped $1$-deck. The same claim holds for the case when one punctures the last bit from both $\bfx^{(1)}$ and $\bfy^{(1)}$. Finally, the claim is true when one punctures both the first and the last bit from both strings. Hence, for $i=1$,
\begin{align}
\cB^{(i)}(\bfx^{(i)}) &= \cB^{(i)}(\bfy^{(i)}),\hspace{2mm} \cB^{(i)}_{LR}(\bfx^{(i)}) = \cB^{(i)}_{LR}(\bfy^{(i)}) \label{eq:equalityB}\\  
\cB^{(i)}_L(\bfx^{(i)}) &= \cB^{(i)}_L(\bfy^{(i)}),\hspace{2mm} \cB^{(i)}_R(\bfx^{(i)}) = \cB^{(i)}_R(\bfy^{(i)}). \nonumber
\end{align}
Let $G^*(k)$ be the smallest integer $n$ such that there exist two distinct binary strings of length $n$, $\bfx^{(k)}$ and $\bfy^{(k)}$, for which~\eqref{eq:equalityB} holds for the case $i=k$.

\begin{theorem} With $\bfx^{(1)}$ and $\bfy^{(1)}$ defined as in (\ref{eq:x1}) and
\begin{align}
\bfx^{(k)} &= \left( 0, \bfx^{(k-1)}, 0, 0, \bfy^{(k-1)}, 0 \right ), \label{eq:rec} \\
\bfy^{(k)} &= \left( 0, \bfy^{(k-1)}, 0, 0, \bfx^{(k-1)}, 0 \right ), \nonumber
\end{align}
defined recursively, we have that \eqref{eq:equalityB} holds for all i. As a result, $G^*(k) \leq 4(2^k-1)$ and $G(k) \leq 4(2^k-1)-2$.
\end{theorem}
\begin{IEEEproof} We split the proof into four subproofs, in order to show that each of the four conditions in~\eqref{eq:equalityB} hold for $i=k$ if they hold for $i=k-1$. We do this by partitioning each deck in~\eqref{eq:equalityB} with respect to whether each padded $0$ is included in a subsequence or not, and by showing that there exists a correspondence between each pair of decks. The bound follows since the length of $\bfx^{(k)}$ equals $4(2^k-1)$  and $G(k) \leq G^*(k)-2$, given that one can remove the padded $0$s.

\subsection*{Part 1: Proof that $\cB_{LR}^{(k)}\left( \bfx^{(k)} \right) = \cB_{LR}^{(k)}\left( \bfy^{(k)} \right)$.}
By definitions \eqref{eq:defs3} and \eqref{eq:rec}, this is equivalent to showing that
\begin{align*}
\cB^{(k)}\left(\bfx^{(k-1)}, 0, 0, \bfy^{(k-1)} \right) = \cB^{(k)}\left(\bfy^{(k-1)}, 0, 0, \bfx^{(k-1)} \right).
\end{align*}

We can partition $\cB^{(k)}\left(\bfx^{(k-1)}, 0, 0, \bfy^{(k-1)} \right)$ depending on which of the two $0$s, if any, is included in the subsequence:
\begin{enumerate}
\item $\{{\left( \cD^{(k_1)}\left(\bfx^{(k-1)}\right), \cD^{(K-k_1)}\left(\bfy^{(k-1)}\right)  \right)\}}$, $K \leq k$;
\item $\{\left( \cD_R^{(k_1)}\left(\bfx^{(k-1)}\right), 0,  \cD^{(K-k_1)}\left(\bfy^{(k-1)}\right) \right)\}$, $K \leq k-1$;
\item $\{\left( \cD^{(k_1)}\left(\bfx^{(k-1)}\right), 0,  \cD_L^{(K-k_1)}\left(\bfy^{(k-1)}\right) \right)\}$, $K \leq k-1$,
\end{enumerate}
where $k_1$ varies from $0$ to $K$, for all $K$. First, we consider the case where neither of the two $0$s is used and show that
\begin{align}\label{eq:show1}&\{\left( \cD^{(k_1)}\left(\bfy^{(k-1)}\right), \cD^{(K-k_1)}\left(\bfx^{(k-1)}\right)  \right)\}=\\\notag &\{\left( \cD^{(k_1)}\left(\bfx^{(k-1)}\right), \cD^{(K-k_1)}\left(\bfy^{(k-1)}\right)  \right)\},
\end{align}
for any $K \leq k$. In this case, each string comprises $k_1 \leq K$ symbols from $\bfx^{(k-1)}$ and $K-k_1$ symbols from $\bfy^{(k-1)},$ where $K \leq k$ denotes the length of the resulting string. When $k_1=K$ or $k_1=0$, (\ref{eq:show1}) holds, since the subsequences $\bfx^{(k-1)}, \bfy^{(k-1)}$ appear in both $\bfx^{(k)}$ and $\bfy^{(k)}$. Otherwise, when $0 < k_1 < K$, since the subsequences $\bfx^{(k-1)}, \bfy^{(k-1)}$ appear (and are ``nonadjacent") in both $\bfx^{(k)}$ and $\bfy^{(k)}$, $\cD^{(k_1)}\left(\bfx^{(k-1)} \right) = \cD^{(k_1)}\left(\bfy^{(k-1)} \right)$, and $\cD^{(K-k_1)}\left(\bfy^{(k-1)} \right) = \cD^{(K-k_1)}\left(\bfx^{(k-1)} \right)$ (since both $k_1 < k$ and $K-k_1 < k$), it follows that (\ref{eq:show1}) also holds for $0 < k_1 < K$.
  
The multiset of subsequences covered by case 2 contains strings that are formed by concatenating $k_1$ bits from $\bfx^{(k-1)}$, the first $0$ between the subsequences $\bfx^{(k-1)}$ and $\bfy^{(k-1)}$ and $K-k_1$ bits from $\bfy^{(k-1)}$. Next, we show that 
 \begin{align}\label{eq:show2}
&\{\left( \cD_R^{(k_1)}\left(\bfx^{(k-1)}\right), 0,  \cD^{(K-k_1)}\left(\bfy^{(k-1)}\right) \right)\}= \\ \notag
&\{\left( \cD_R^{(k_1)}\left(\bfy^{(k-1)}\right), 0,  \cD^{(K-k_1)}\left(\bfx^{(k-1)}\right) \right)\}
 \end{align}
 for $K \leq k-1$. Since $K \leq k-1$, we have $\cD_R^{(k_1)}(\bfx^{(k-1)}) = \cD_R^{(k_1)}(\bfy^{(k-1)})$ and $\cD^{(K-k_1)}(\bfy^{(k-1)}) = \cD^{(K-k_1)}(\bfx^{(k-1)})$, which implies that we can form strings by concatenating $k_1$ bits from $\bfy^{(k-1)}$, the first $0$ between the substrings $\bfy^{(k-1)}$ and $\bfx^{(k-1)}$, and $K-k_1$ bits from $\bfx^{(k-1)}$. Thus, (\ref{eq:show2}) holds. 
 
 Using the same approach, it can be shown that 
\begin{align}\label{eq:toshow3}
&\{\left( \cD^{(k_1)}\left(\bfx^{(k-1)}\right), 0,  \cD_L^{(K-k_1)}\left(\bfy^{(k-1)}\right) \right)\} =\\ \notag
&\{\left( \cD^{(k_1)}\left(\bfy^{(k-1)}\right), 0,  \cD_L^{(K-k_1)}\left(\bfx^{(k-1)}\right) \right)\},
 \end{align}
for any $K \leq k-1$. From (\ref{eq:show1}), (\ref{eq:show2}), and (\ref{eq:toshow3}), it then follows that $\cB_{LR}^{(k)}\left( \bfx^{(k)} \right)=\cB_{LR}^{(k)}\left( \bfy^{(k)} \right)$.

\subsection*{Part 2: Proof that $\cB_{L}^{(k)}\left( \bfx^{(k)} \right) = \cB_{L}^{(k)}\left( \bfy^{(k)} \right)$.} 
By definitions \eqref{eq:defs1} and \eqref{eq:rec}, this is equivalent to showing that

\begin{align*}
\cB^{(k)}\left(\bfx^{(k-1)}, 0, 0, \bfy^{(k-1)},0 \right) = \cB^{(k)}\left(\bfy^{(k-1)}, 0, 0, \bfx^{(k-1)},0 \right).
\end{align*}

We first partition $\cB_{L}^{(k)}\left( \bfx^{(k)} \right)$ into two multisets: The first contains subsequences that include the last (trailing) $0$ while the second contains those which do not (equivalent to $\cB_{LR}^{(k)}(\bfx^{(k)})$). The first multiset can be partitioned into three classes:
\begin{enumerate}[leftmargin=0.6cm]
\item $\{\left( \cD^{(k_1)}\left(\bfx^{(k-1)}\right), \cD_R^{(K-k_1)}\left(\bfy^{(k-1)}\right)  , 0 \right)\}$, $K \leq k-1$;
\item $\{\left( \cD_R^{(k_1)}\left(\bfx^{(k-1)}\right), 0,  \cD_R^{(K-k_1)}\left(\bfy^{(k-1)}\right), 0 \right)\}$, $K \leq k-2$;
\item $\{\left( \cD^{(k_1)}\left(\bfx^{(k-1)}\right), 0,  \cD_{LR}^{(K-k_1)}\left(\bfy^{(k-1)}\right), 0 \right)\}$, $K \leq k-2$.
\end{enumerate}
Using an almost identical argument as the one described in \textit{Part 1}, one can show that $\cB_{L}^{(k)}\left( \bfx^{(k)} \right) = \cB_{L}^{(k)}\left( \bfy^{(k)} \right)$.

\subsection*{Part 3: Proof that $\cB_{R}^{(k)}\left( \bfx^{(k)} \right) = \cB_{R}^{(k)}\left( \bfy^{(k)} \right)$.} 
The proof of this case follows by symmetry from \textit{Part 2}.
\subsection*{Part 4: Proof that $\cB^{(k)}\left( \bfx^{(k)} \right) = \cB^{(k)}\left( \bfy^{(k)} \right)$.} The final step in the proof is to show that
\begin{align*}
\cB^{(k)}\left(0, \bfx^{(k-1)}, 0, 0, \bfy^{(k-1)}, 0 \right) \\
= \cB^{(k)}\left(0, \bfy^{(k-1)}, 0, 0, \bfx^{(k-1)}, 0 \right)
\end{align*}
Using a similar approach as before, we now partition the subsequences in $\cB^{(k)}(\bfx^{(k)})$ according to whether they
\begin{enumerate}
    \item contain the leading $0$, but not the trailing $0$;
    \item contain the trailing $0$, but not the leading $0$;
    \item contain neither the trailing nor the leading $0$;
    \item contain both the leading and the trailing $0$.
\end{enumerate}
This is equivalent to:
\begin{enumerate}
    \item $\cB^{(k)}_{R}\left( \bfx^{(k)} \right)\setminus\cB^{(k)}_{LR}\left( \bfx^{(k)} \right)$;
    \item $\cB^{(k)}_{L}\left( \bfx^{(k)} \right)\setminus\cB^{(k)}_{LR}\left( \bfx^{(k)} \right)$;
    \item $\cB^{(k)}_{LR}\left( \bfx^{(k)} \right)$;
    \item $\cB^{(k)}\left( \bfx^{(k)} \right)\setminus(\cB^{(k)}_{R}\left( \bfx^{(k)} \right)\cup(\cB^{(k)}_{L}\left( \bfx^{(k)} \right))$.
\end{enumerate}
From the first three parts of the proof, we know that the first three multisets are the same for $\bfx^{(k)}$ and $\bfy^{(k)}$. We only need to prove that the fourth multiset is the same as well. Again we partition the multiset of interest into three classes:
\begin{enumerate}
\item $\{\left(0, \cD_{L}^{(k_1)}\left(\bfx^{(k-1)}\right), \cD_R^{(K-k_1)}\left(\bfy^{(k-1)}\right), 0 \right)\}$;\\
\item $\{\left(0, \cD_{LR}^{(k_1)}\left(\bfx^{(k-1)}\right), 0,  \cD_R^{(K-k_1)}\left(\bfy^{(k-1)}\right), 0 \right)\}$;\\
\item $\{\left(0, \cD_{L}^{(k_1)}\left(\bfx^{(k-1)}\right), 0,  \cD_{LR}^{(K-k_1)}\left(\bfy^{(k-1)}\right), 0 \right)\}$,
\end{enumerate}
where for case 1, $K \leq k-2$, and for cases 2 and 3, $K \leq k-3$. Using similar arguments as before completes the proof.
\end{IEEEproof}
Using a similar approach, we can extend the bound to the $s$-gapped case to get $G_s(k)\leq (5s-2)2^{k-1}-5s+4$. This is done by adding $s-1$ $0$s on the outside and $s$ $0$s between $xy$ and $yx$. We remove $2s-2$ $0$s for the bound since the $s$-gapped $k$-deck does not need to satisfy the extra conditions required by the recursive construction.

We numerically computed $G(k)$ for $k=2,3,4$. The results are displayed below,
\begin{table}[h]
\centering
\begin{tabular} {| c | c | c |} 
 \hline
  $k$ & $G(k)$ &  Confusable pairs (examples)\\ 
 \hline
  2 & 6 & (0,1,0,0,1,1), (0,0,1,1,0,1) \\ 
 \hline
  3 & 13 &(1,1,0,1,1,1,1,0,1,0,1,1,1), (1,1,1,0,1,0,1,1,1,1,0,1,1)\\ 
 \hline
  4 & 24 & (1,1,0,0,1,1,0,1,0,1,0,1,0,0,1,1,0,0,1,1,0,1,0,0),  \\ 
   &  & (1,1,0,1,0,0,1,1,0,0,1,1,0,1,0,1,0,1,0,0,1,1,0,0) \\
 \hline
\end{tabular} \label{tab:values}
\end{table}
which clearly indicate that the upper bound $4(2^k-1)-2$ is loose for larger values of $k$: For $k=4$, the bound equals $58$ while the correct value is only $24$. Also, the exact values of $G(k)$ are significantly larger than those for the ungapped case, for which we know that $S(k)=4,7,12$ (compared to $G(k)=6,13,24$) for $k=2,3,4$, respectively. We therefore turn our attention to improving the bound on $G(k)$ using more sophisticated counting arguments.

\section{Improved Upper Bounds for Gapped $k$-Decks} \label{sec:improve}

We find the following definitions and notation from~\cite{Dud} useful for our subsequent derivations.
Let $\Gamma=\{X,Y\}$ and let $J$ denote a ``wildcard''. For integers $0\leq r\leq k$ let
\vspace{-4mm}
\begin{equation*}
    U_r(k)=\{\bfw \in\bigcup_{j=r}^{k} (\Gamma\cup\{J\})^j:\; \bfw \text{ has exactly $r$ non-$J$ symbols}\}.
\end{equation*}
For $t\geq1$ and $k_1\geq\cdots\geq k_t\geq t$, let
\begin{equation}
    U(k_1,\ldots,k_t)=U_1(k_1)\cup U_2(k_2)\cup\cdots\cup U_t(k_t).
\end{equation}
We restrict our attention to $U(k_1,k_2)=U_1(k_1)\cup U_2(k_2),$ the set of all strings of length at most $k_1$ that have exactly one non-$J$ character and the set of strings of length at most $k_2$ that have exactly two non-$J$ characters.

When we refer to the multiplicity with which a string $\bfw$ that contains wildcards ($J$'s) occurs as a subsequence of a string $\bfp$ that contains no wildcards (denoted by $N(\bfw,\bfp)$), we  map each wildcard to either $X$ or $Y$. For example, if $\bfw=(J,X)$ and $\bfp=(Y,X,Y,X)$, we have $N(\bfw,\bfp)=4$ because $(X,X)$ and $(Y,X)$ occur as subsequences of $\bfp$ with multiplicity $1$ and $3$, respectively. Let $\bfp$ and $\bfq$ be two binary strings. We write $\bfp\sim^{U_r(k)}\bfq$ if $N(\bfw,\bfp)=N(\bfw,\bfq)$ for all $\bfw\in U_r(k)$. In addition, we write $\bfp\sim^{U(k_1,k_2)}\bfq$ if $N(\bfw,\bfp)=N(\bfw,\bfq)$ for all $\bfw\in U(k_1,k_2)$.

Next, let $S_{U}(k_1)$ be the smallest integer $m$ for which there exist distinct strings $\bfp$ and $\bfq$ of length $m$ such that $\bfp\sim^{U_r(k_1)}\bfq$. Similarly, let $S_{U}(k_1,k_2)$ be the
smallest integer $m$ for which there exist distinct strings $\bfp$ and $\bfq$ of length $m$ such that $\bfp\sim^{U(k_1,k_2)}\bfq$. The following lemma is used in our subsequent derivations. 
\begin{lemma}\cite{Dud}\label{Dudik}
Let $k_1\geq k_2\geq 2$ and $\kappa=k_1^2+k_2^2(k_2-1)/2$. Then $S_{U}(k_1,k_2)\leq \kappa (\lg \kappa + \lg\lg \kappa+1) = (1+o(1)) \kappa \lg \kappa$.
\end{lemma}

Let $N_g(\bfw,\bfx)$ be the number of times a string $\bfw$ appears as a gapped subsequence of $\bfx$ (i.e., so that all indices in $\bfx$ are nonadjacent). When $N_g(\bfw,\bfx)=N_g(\bfw,\bfy)$ for all strings $\bfw$ of length $\leq k$, then we write $\bfx\sim^{k_{(g)}}y$, i.e. $\cB^{(k)}(\bfx)=\cB^{(k)}(\bfy)$. 

Let $\Gamma=\{X,Y\}$ and let $\Sigma$ be an arbitrary alphabet. For a finite-length string $\bfx$ over $\Sigma$, define $\bfx_0$ to be the string obtained by padding $x$ with one $0$ at both ends. For a finite-length string $\bfp$ over $\Gamma$ and two finite-length strings $\bfx$ and $\bfy$ over $\Sigma$, let $h_{\bfx,\bfy}(\bfp)$ be the string obtained from $\bfp$ by replacing each $X$ by the string $\bfx_0$ and each $Y$ by the string $\bfy_0$. For example, If $\bfp=(X,Y)$ and $\bfx=(0,1,0,1)$ and $\bfy=(1,1,0,0)$, then $\bfx_0=(0,0,1,0,1,0)$, $\bfy_0=(0,1,1,0,0,0)$ and $h_{\bfx,\bfy}(\bfp)=(0,0,1,0,1,0,0,1,1,0,0,0)$. We are now ready to prove an analogue of Lemma 9 from~\cite{Dud} for the case of gapped $k$-decks.\vspace{-0.03in}

\begin{lemma} \label{lem}
Let $\bfx$ and $\bfy$ be two distinct strings in $\Sigma^n$ such that $\bfx\sim^{k_{(g)}}\bfy$, $\cB^{(k)}_L(\bfx)=\cB^{(k)}_L(\bfy)$, $\cB^{(k)}_R(\bfx)=\cB^{(k)}_R(\bfy)$ and $\cB^{(k)}_{LR}(\bfx)=\cB^{(k)}_{LR}(\bfy)$. Let $\bfp$ and $\bfq$ be two distinct binary strings in $\Gamma^m$, such that for some $\sigma\in\{0,1,2\},$ we have $\bfp\sim^{U(2k+\sigma,k+\sigma)}\bfq$. Then, $h_{\bfx,\bfy}(\bfp)$ and $h_{\bfx,\bfy}(\bfq)$ are distinct and we have $h_{\bfx,\bfy}(\bfp)\sim^{3k+\sigma_{(g)}}h_{\bfx,\bfy}(\bfq)$. The same result holds when puncturing $h_{\bfx,\bfy}(\bfp)$ and $h_{\bfx,\bfy}(\bfq)$ on the left, right, and on both sides by one bit.
\end{lemma}\vspace{-0.03in}
\begin{IEEEproof}
Due to space limitations, we only provide a sketch of the proof. Let $\bfw$ be a string of length at most $3k+\sigma$ in $\Sigma$, and $p=(p_1,\ldots,p_m)$. The idea of the original proof \cite{Dud} for the ungapped case is as follows: Each mapping that takes $\bfw$ to $h_{\bfx,\bfy}(\bfp)$ (as a subsequence) defines a splicing of $\bfw$ of the form $\bfw=(\bfw_1,\bfw_2,\ldots,\bfw_m)$, where $\bfw_i$ is the preimage of $h_{\bfx,\bfy}(\bfp_i)$. Note that some $\bfw_i$ strings may be empty. Hence, we can write the set of all mappings which take $\bfw$ to $h_{\bfx,\bfy}(\bfp)$ as the union of direct products (see~\cite{Dud} for the specific notation) $\bigcup_{t\geq1}\bigcup_{l}\bigcup_r\prod_{i=1}^t \cN(\bfw_{l,i},h_{\bfx,\bfy}(\bfp_{r(i)}))$, where $t$ is the number of nonempty segments, the second union is taken over all functions $l$ which partition $\bfw$ into $t$ nonempty segments for a fixed $t$, the third union is taken over all functions $r$ mapping the $t$ chosen nonempty segments in $\bfw$ to $t$ of the segments in $h_{\bfx,\bfy}(\bfp)$ (equivalently, each $r$ corresponds to a way in which we pick $t$ out of $m$ segments of $h_{\bfx,\bfy}(\bfp)$), and finally, $\cN(\bfw_{l,i},h_{\bfx,\bfy}(\bfp_{r(i)}))$ denotes the set of mappings taking the $i$-th nonempty segment of $\bfw$ into the corresponding chosen segment in $h_{\bfx,\bfy}(\bfp)$. Note that for every mapping $f$ which takes $\bfw$ to $h_{\bfx,\bfy}(\bfp)$, there is a specific $t$, $l$ and $r$ that correspond to $f$. Furthermore, $f$ is the direct product of $t$ mappings, each taking one of the $t$ nonempty segments of $\bfw$ to the corresponding segment of $h_{\bfx,\bfy}(\bfp)$ (which is all uniquely determined by fixing $t$, $l$ and $r$). Now, by converting this expression into a corresponding sum, we get $N(\bfw,h_{\bfx,\bfy}(\bfp))=\sum_{t\geq1}\sum_{l}\sum_r\prod_{i=1}^t N(\bfw_{l,i},h_{\bfx,\bfy}(\bfp_{r(i)}))$. Since $\bfw$ has length at most $3k+\sigma\leq3k+2$, it has at most two segments of length $>k$. Therefore, we have three types of $l$'s (i.e., ways of partitioning $\bfw$ into $t$ nonempty segments for a fixed $t$): The ones with no segments of length $>k$, the ones with one such segment and the ones with two such segments. This means that $N(\bfw,h_{\bfx,\bfy}(\bfp))$ is a triple sum. Since $\bfp\sim^{U(2k+\sigma,k+\sigma)}\bfq$, after some calculations we get $N(\bfw,h_{\bfx,\bfy}(\bfp))=N(\bfw,h_{\bfx,\bfy}(\bfq))$. 

To adapt this procedure for the gapped case, we need to show that $N_g(\bfw,h_{\bfx,\bfy}(\bfp))=N_g(\bfw,h_{\bfx,\bfy}(\bfq))$. The first difference is that each gapped mapping that takes $\bfw$ to $h_{\bfx,\bfy}(\bfp)$ (as a gapped subsequence) defines a splicing of $\bfw$ of the form $\bfw=\bfw_1\bfz_1\bfw_2\bfz_2\ldots \bfz_{m-1}\bfw_m$, where $\bfw_i$ is the preimage of $h_{\bfx,\bfy}(\bfp_i)$  and $\bfz_i$ is the preimage of the $i$-th pair of $0$s that we added between the $\bfx$ and between the $\bfy$ strings and between the $\bfx$ and $\bfy$ strings when constructing $h_{\bfx,\bfy}(\bfp)$. Again, note that some $\bfw_i$ and $\bfz_i$ strings may be empty. 

However, an important difference between the gapped and ungapped case is that we need to consider different cases based on whether $\bfz_i$ is empty or not, for all indices $i$. This is because if $\bfz_1$ is nonempty, for example, then we need to make sure that we do not use the rightmost bit in $h_{\bfx,\bfy}(\bfp_1)$ (or the leftmost bit in $h_{\bfx,\bfy}(\bfp_2)$, depending on whether the gapped mapping takes $\bfz_1$ to the first or second $0$ from the pair of $0$s in $h_{\bfx,\bfy}(\bfp)$ that are positioned between $h_{\bfx,\bfy}(\bfp_1)$ and $h_{\bfx,\bfy}(\bfp_2)$). This case ($\bfz_1$ nonempty) gives rise to several additional cases that need to be considered, depending on which of the strings $\bfz_2,\bfz_3,\ldots,\bfz_{m-1}$ are empty. In other words, we can write out $h_{\bfx,\bfy}(\bfp)=(0,s_1,0,0,s_2,0,0,...,0,0,s_m,0)$ where each $s_i\in \{\bfx,\bfy\}$. Any gapped subsequence will then be of the form $(J_0,\cD_{\alpha_1}^{(i_1)}(s_1),J_1,\cD_{\alpha_2}^{(i_2)}(s_2),J_2,...,J_{m-1},\cD_{\alpha_m}^{(i_m)}(s_m),J_m)$ where $J_j\in\{0,\emptyset\}$ and $\alpha_j\in \{\emptyset, L, R, LR\}$. Here, each $J_j$ represents a $0$ in the padding and whether it is a part of the subsequence or not. Then, depending on whether or not we use the padding, we puncture $s_j$ on the left, right, both, or neither. This is captured by the indices $\alpha_j$'s. We also have that $\sum_{j=1}^m i_j+\sum_{j=1}^m |J_j|\leq3k+\sigma$. Since by our assumptions $\cD_{\alpha}^{(i)}(p_j) = \cD_{\alpha}^{(i)}(q_j)$ for all $i\leq k$, $1\leq j \leq m$, $\alpha_j\in \{\emptyset, L, R, LR\}$, we have an equivalence between $h_{\bfx,\bfy}(\bfp)$ and $h_{\bfx,\bfy}(\bfq)$ for each $i_j\leq k$. By the summation constraint, there are at most two indices $j$ such that $i_j\geq k+1$. Let us consider the case when there is exactly one such $j$, denoted by $j^*$. In this case we have to pick fewer than $2k+\sigma$ of the remaining characters to obtain the final subsequence. We can also divide $h_{\bfx,\bfy}(\bfp)$ into subblocks of the form $(0,s_j,0)$, i.e., we can splice $\bfw$ into a collection of $\bfw_j$'s, where $\bfw_j$ is the string mapped to one of the blocks and $|\bfw_j|<k$ for $j\neq j^*$. The multiplicity of $\bfw$ can be seen to be $N_g(\bfw_{j^*},(s_{j^*})_{\alpha_{j^*}})\prod_{j\neq j^*} N_g(\bfw_j,(\bfx)_{\alpha_j})$, where $\alpha$ once again depends on whether the bit used for padding is included in the subsequence. By our assumption we have $N(J^{a_1} A J^{a_2}, \bfp)=N(J^{a_1} A J^{a_2}, \bfq),$ where $J^{a_1}$ is a sequence of $a_1$ concatenated wildcard characters, $A\in \{x,y\}$, $a_1+a_2<2k+\sigma$. Hence, there are equally many $s_{j^*}$'s in $h_{\bfx,\bfy}(\bfp)$ and $h_{\bfx,\bfy}(\bfq)$. Using similar arguments and the fact that $N(J^{a_1} A J^{a_2} B J^{a_3}, \bfp)=N(J^{a_1} A J^{a_2} B J^{a_3}, \bfq)$ we can also prove the equivalence for the case of two indices $j$ for which $i_j\geq k+1$. This leads to $N_g(\bfw,h_{\bfx,\bfy}(\bfp))=N_g(\bfw,h_{\bfx,\bfy}(\bfq))$.
\end{IEEEproof}
The lemma gives rise to the following important Corollary.
\begin{corollary}
For every $\sigma\in\{0,1,2\}$, one has $G^*(3k+\sigma)\leq (G^*(k)+2)(S_{U}(2k+\sigma,k+\sigma))$. \label{inequalityG}
\end{corollary}
Combining the above corollary with Lemma~\ref{Dudik} and Lemma~\ref{lem} leads to an upper bound for $G(k)$ as follows. First, we set 
$$\kappa=(2k+\sigma)^2+(k+\sigma)^2(k+\sigma-1)/2,$$
which equals
\begin{align*}
&\frac{1}{2}(\sigma^3+(3k+1)\sigma^2+3\sigma(k^2+2k)+(1+7/k)k^3)\\
=&(\frac{1}{2}+o(1))\,k^3,
\end{align*}
in Lemma~\ref{Dudik} to obtain
\begin{equation}
    S_{U}(2k+\sigma,k+\sigma)\leq C(k) k^3 log_3 k,
    \label{S_Ubound}
\end{equation}
where $C(k) = \frac{3\lg3}{2}+o(1)$. We also have that $C(k)\leq 10$ for $k\geq 9$, $C(k)\leq 3$ for $k\geq 3^5$~\cite{Dud}. Using the inequality from Corollary \ref{inequalityG} and~\eqref{S_Ubound}, we set $k_{0}=k$ and for $i>0$, $k_i =\lfloor k_{i-1}/3 \rfloor\leq k/3^i$. We stop the recursion with $i=i_0$, where $k_{i_0}\leq 4$ (so $G^*(k_{i_0})\leq 4(2^4-1) = 60$) and get
\vspace{-2mm}
\begin{align}
    G^*(k) &\leq G^*(k_{i_0})\prod_{i=1}^{i_0}S_U(2k_i+\sigma_i,k_i+\sigma_i)\\ \notag
    &+2\sum_{i=1}^{i_0}\prod_{j=1}^{i}S_U(2k_i+\sigma_i,k_i+\sigma_i).
\end{align}
Combining the above bound with that on $S_U$ we obtain
\vspace{-2mm}
\begin{align*}
G^*(k) &\leq G^*(k_{i_0})\prod_{i=1}^{i_0}C(k_i) (k_i)^3 log_3 (k_i)\\ &+2\sum_{i=1}^{i_0}\prod_{j=1}^{i}C(k_i) (k_i)^3 log_3 (k_i)\\
&\leq 3^{\log_3(60) + \sum_{i=1}^{\lceil\log_3(k/4)\rceil}\left[O(1)+3(\log_3 k-i)+\log_3(\log_3 k-i)\right]}\\
&+\sum_{i=1}^{\lceil\log_3(k/4)\rceil}3^{\log_3 2+\sum_{j=1}^i\left[O(1)+3(\log_3 k-j)+\log_3(\log_3 k-j)\right]}\\
&=3^{O(1) + O(\log_3 k)+O(\log_3^2 k)+O(\log_3 k\log_3 \log_3 k)}\\
&+O(\log_3(k))3^{O(1) + O(\log_3 k)+O(\log_3^2 k)+O(\log_3 k\log_3 \log_3 k)}\\
&=O(\log_3(k))3^{O(\log_3^2 k)}.
\end{align*}
\vspace{-6mm}

Since we have $G(k)\leq G^*(k)-2$, we also have $G(k)\leq O(\log_3(k))3^{O(\log_3^2 k)}$. By bounding $G^*(k)$ we also obtain
\vspace{-1mm}
\begin{align*}
    G(k) &\leq (4(2^{k/3}-1)+2)C(k/3)(k/3)^3\log_3(k/3) - 2\\
    &\leq 4/27*2^{k/3}*k^3\log_3(k/3)*C(k/3)-2
\end{align*}
for $k\geq 28$. Since in this case $C(k/3)\leq 10$, we arrive at
\vspace{-1mm}
\begin{align}
    G(k)\leq 1.482*1.26^k*k^3\log_3(k/3)-2,
\end{align}

\vspace{-2mm}

In comparison, the general bound for the ungapped case, derived in~\cite{Dud}, reads as
\vspace{-1mm}
$$S(k) \leq 1.2 \Gamma(\log_3\,k)\,3^{(3/2)\,\log_3^2\,k-(1/2)\log_3\,k},\; k\geq 85.$$
The bounds are summarized in the tables below.
\begin{table}[h]
\centering
\begin{tabular} {|| c || c || c ||} 
 \hline
  $k$ & Bound\\ 
 \hline
  2-4 &  Exact values: 6,13,24 \\ 
 \hline
  5-27 & $4(2^{k}-1)$\\ 
 \hline
  $>27$ & $1.482*1.26^k*k^3\log_3(k/3)-2$\\ 
 \hline
\end{tabular}
\vspace{-1.8em}
\end{table}
\begin{table}[ht]
\begin{adjustbox}{width=\columnwidth,center}
\begin{tabular} {|| c || c | c | c | c | c | c  ||} 
 \hline
  $k$ & 28 & 29 & 30 & 31 & 32 & 33\\ 
 \hline
  $G(k)\leq$ & 42742211&60773950&86039831&121319982&170424514&238563374 \\
 \hline
\end{tabular}
\end{adjustbox}
\end{table}


\end{document}